\documentclass[11pt]{article}
\usepackage[T1]{fontenc}
\usepackage[a4paper]{anysize}\marginsize{2.0cm}{2.0cm}{1.0cm}{1.0cm}
\pdfpagewidth=\paperwidth \pdfpageheight=\paperheight
\usepackage{pgf,tikz,hyperref}
\usetikzlibrary{arrows, arrows.meta}
\usepackage{amsfonts,amssymb,amsthm,amsmath,eucal}
\usepackage{graphicx}
\usepackage{ltablex}
\usepackage{caption, booktabs}
\usepackage{longtable}
\usepackage{enumitem}
\usepackage{multirow}
\usepackage{makecell}
\usepackage{pgfplots}

\usepackage{authblk}

\pgfplotsset{compat=1.15}
\usepackage{mathrsfs}

\theoremstyle{plain}
\newtheorem{thm}{Theorem}[section]

\newtheorem{lemma}[thm]{Lemma}
\newtheorem{proposition}[thm]{Proposition}

\newtheorem{thmABC}{Theorem}

\theoremstyle{definition}
\newtheorem{definition}[thm]{Definition}
\newtheorem{remark}[thm]{Remark}
\newtheorem{example}[thm]{Example}

\newtheorem{fact}[thm]{Fact}

\newtheorem{thevarthm}[thm]{\varthmname}

\newenvironment{varthm*}[1]{\trivlist\item[]{\bf #1.}\it}{\endtrivlist}

\newcommand\be{\begin{eqnarray*}}
\newcommand\ee{\end{eqnarray*}}

\renewcommand\P{\mathbb P}

\newcommand\newop[2]{\def#1{\mathop{\rm #2}\nolimits}}
\newop\edim{edim}
\newop\Zeroes{Zeroes}
\newop\Jac{Jac}
\newop\Ass{Ass}
\newop\SL{SL}
\newop\PGL{{\P}GL}
\newop\Km{Km}
\newop\reg{reg}

\newcolumntype{L}{>{$}l<{$}}

\newcommand\keywords[1]{{\renewcommand\thefootnote{}\footnotetext{\textit{Keywords:} #1.}}}
\newcommand\subclass[1]{{\renewcommand\thefootnote{}\footnotetext{\textit{Mathematics Subject Classification (2020):} #1.}}}

\newcommand{\n}{\boldsymbol{\mathcal{N}}}
\renewcommand{\t}{\boldsymbol{\mathcal{T}}}
\renewcommand{\tt}{\boldsymbol{\mathcal{T}^2}}
\renewcommand{\a}[1]{\boldsymbol{\mathcal{A}^{#1}}}

\begin{document}

\author[1]{\L{}ukasz Merta}
\author[1]{Filip Zieli\'nski}
\author[1]{Marcin Zieli\'nski}
\title{On free arrangements of three conics}
\affil[1]{University of the National Education Commission, Krakow}
\date{\today}
\maketitle
\thispagestyle{empty}
\begin{abstract}
    We give a complete classification of free arrangement of three smooth conics on complex projective plane  admitting only ${\rm ADE}$ singularities and $J_{2,0}$ singularities.
\end{abstract}
\subclass{MSC 14C20 \and MSC 14N20}

\section{Introduction}
The aim of this paper is to classify all free arrangements (up to the projective equivalence) of three smooth conics in $\mathbb{P}^2_{\mathbb{C}}$ admitting only certain quasi-homogeneous singularities. Let us recall that an arrangement of plane curves is called {\it free} if its associated module of derivations is a free module over the coordinate ring of the complex projective plane. The free arrangements of conics with quasi-homogeneous singularities were studied before by some authors, see for example \cite{DJP, PP23, Pokora2023, SCHENCK}, but still our knowledge and understanding of these arrangements is far from being complete, even in a very natural setting of conic arrangements with only ${\rm ADE}$ singularities. Free curves are not common objects because they are generally difficult to find or construct. 

The inspiration for our study of arrangements of three smooth conics is \cite[Theorem 2.2]{Pokora2023} which states that if an arrangement of $k \geq 2$ smooth conics with some simple singularities is free, then $k \in \{2,3,4\}$. This restriction is, to some extent, surprising and it is very natural to ask whether we can completely classify such arrangements.

It is known that there exists only one (up to the projective equivalence) $1$-parameter family of two smooth conics that is free \cite[Proposition 5.5]{DIMCA202377}, i.e., every member of this family is free. Additionally, two free arrangements of three smooth conics are also known. One of them is the arrangement with the defining equations
\begin{equation*}
\begin{array}{l}
Q_{1} \colon X^{2}+Y^{2}-Z^{2} = 0, \\
Q_{2} \colon 2X^2 + Y^2 +2XZ = 0, \\
Q_{3} \colon 2X^2 + Y^2 - 2XZ = 0.
\end{array}
\end{equation*}
This arrangement is well-known in algebraic geometry and it is called Persson's triconical arrangement \cite{Persson1982}. 
The other one is the arrangement presented in \cite{Pokora2023}, with three conics defined by the following equations:
\begin{equation*}
\begin{array}{l}
Q_{1} \colon -3X^2 + XY + YZ + ZX = 0, \\
Q_{2} \colon -3Y^2 + XY + YZ + ZX = 0, \\
Q_{3} \colon -3Z^2 + XY + YZ + ZX = 0.
\end{array}
\end{equation*}
We are also aware of examples of free arrangements of four smooth conics, for instance one based on Persson's triconical arrangement (see \cite[Example 3.3]{MAXIMIZING}). However, we are not aware of any free arrangements of four smooth conics admitting only singularities of type $A_1$, $A_3$, $A_5$, $A_7$ and $D_4$ (this notation is according to Arnold's classification of singularities \cite{Arnold1976}).

Our research project focuses on finding all free arrangements of three and four smooth conics with ${\rm ADE}$ singularities, and in this article we make the first step into this direction, namely we deliver a classification result on free arrangements of three smooth conics, see Theorem \ref{thm:main1}. Additionally, we show that there are no free arrangements containing a singularity of type $J_{2,0}$, which is a quasi-homogeneous singularity that is not ${\rm ADE}$, see Theorem \ref{thm:main2}.

\section{Preliminaries }

For a projective situation, with a point $p\in \mathbb{P}^{2}_{\mathbb{C}   }$ and a homogeneous polynomial $f\in \mathbb{C}[X,Y,Z]$, we take local affine coordinates such that $p=[0:0:1]$ and then the dehomogenization of $f$.

\begin{definition}
Let $p$ be an isolated singularity of a polynomial $f\in \mathbb{C}[X,Y]$. Since we can change the local coordinates, let $p=(0,0)$.
The number 
$$\mu_{p}=\dim_\mathbb{C}\left(\mathbb{C}[[X,Y]] /\bigg\langle \frac{\partial f}{\partial X},\frac{\partial f}{\partial Y} \bigg\rangle\right)$$
is called the Milnor number of $f$ at $p$.\\

The number
$$\tau_{p}=\dim_\mathbb{C}\left(\mathbb{C}[[X,Y]] /\bigg\langle f,\frac{\partial f}{\partial X},\frac{\partial f}{\partial Y}\bigg\rangle \right)$$
is called the Tjurina number of $f$ at $p$.
\end{definition}
Additionally, the total Tjurina number of a given reduced curve $C \subset \mathbb{P}^{2}_{\mathbb{C}}$ is defined as
$$\tau(C) = \sum_{p \in {\rm Sing}(C)} \tau_{p}.$$

Let us denote by $S := \mathbb{C}[X,Y,Z]$ the coordinate ring of  $\mathbb{P}^{2}_{\mathbb{C}}$ and for a homogeneous polynomial $f \in S$ let us denote by $J_{f}$ the Jacobian ideal associated with $f$.
Let $C : f=0$ be a reduced curve in $\mathbb{P}^{2}_{\mathbb{C}}$ of degree $d$ defined by $f \in S$. Denote by $M(f) := S/ J_{f}$ the associated Milnor algebra. 
\begin{definition}
We say that a reduced plane curve $C$ is \emph{an $m$-syzygy curve} when $M(f)$ has the following minimal graded free resolution:
$$0 \rightarrow \bigoplus_{i=1}^{m-2}S(-e_{i}) \rightarrow \bigoplus_{i=1}^{m}S(1-d - d_{i}) \rightarrow S^{3}(1-d)\rightarrow S \rightarrow M(f) \rightarrow 0$$
with $e_{1} \leq e_{2} \leq ... \leq e_{m-2}$ and $1\leq d_{1} \leq ... \leq d_{m}$.
\end{definition}
In the setting of the above definition, the minimal degree of the Jacobian relations among the partial derivatives of $f$ is defined to be ${\rm mdr}(f) := d_{1}$.
\begin{definition}
We say that $C$ is \emph{free} if and only if $m=2$, and then $d_{1}+d_{2}=d-1$.
\end{definition}
A useful criterion for studying the freeness of curves is Du Plessis-Wall theorem \cite{duPlessis}:
\begin{thm}\label{duPles}
    A reduced plane curve $C$ is free if and only if ${\rm mdr}(f) = d_1\leq (d-1)/2$ and
\begin{equation}
(d-1)^{2} - d_{1}(d-d_{1}-1) = \tau(C).
\end{equation}
\end{thm}

We present below the local normal forms of the singularities we are interested in, and these equations are taken from Arnold's paper \cite{Arnold1976}. All of these singularities are quasi-homogeneous.

\begin{table}[h!]
\centering
\begin{tabular}{ll}
$A_{k}$ with $k\geq 1$ & $: \, X^{2}+Y^{k+1}  = 0$, \\
$D_{k}$ with $k\geq 4$ & $: \, Y^{2}X + X^{k-1}  = 0$,\\
$J_{2,0}$ & $: \, X^3 + bX^2Y^2 + Y^6 = 0,\, 4b^3+27\neq 0.$
\end{tabular}
\caption{Local normal forms.}
\label{aloc}
\end{table}

We are going to study arrangements of conics from the perspective of pairs consisting of them. Any arrangement of $k$ smooth conics for $k > 2$ can be split into $\binom{k}{2}$ pairs. In Table \ref{tab: notation}, we show all possible configurations of a pair of smooth conics, along with the graph notation captured from \cite[Table 2.1]{PHD}. In the last column, we introduce our own symbolic notation for each pair. We are going to use this notation throughout the paper.

\begin{table}[h]
\centering
\begin{tabular}{|c|c|c|c|}
    \hline
    Graph notation & Picture & Singularities & Symbolic notation \\
    \hline
    
    \multirow{4}{*}{
        \makecell[c]{ \\
            \begin{tikzpicture}
                \draw [fill] (0,0) circle (3pt);
                \node at (0,-0.5) {$Q_1$};
                \draw [fill] (2,0) circle (3pt);
                \node at (2,-0.5) {$Q_2$};
            \end{tikzpicture}
        }
    } 
    & 
    \multirow{4}{*}{
        \makecell[c]{
            \begin{tikzpicture}[scale=0.7]
                \draw [line width=1.2pt, rotate=30] (0,0) ellipse (2cm and 1cm);
                \draw [line width=1.2pt, rotate=-30] (0,0) ellipse (2cm and 1cm);
            \end{tikzpicture}
        }
    } 
    & & \multirow{4}{*}{$\n$} \\
    & & four $A_1$ & \\ 
    & & (four nodes) & \\ 
    & & & \\ 
    \hline 
    \multirow{4}{*}{
        \makecell[c]{   \\
            \begin{tikzpicture}
                \draw (0,0) -- (2,0);
                \draw [fill] (0,0) circle (3pt);
                \node at (0,-0.5) {$Q_1$};
                \draw [fill] (2,0) circle (3pt);
                \node at (2,-0.5) {$Q_2$};
            \end{tikzpicture}
        }
    }
    & 
    \multirow{4}{*}{
        \makecell[c]{
            \begin{tikzpicture}[scale=0.7]
                \draw [line width=1.2pt, rotate=0] (0,0) ellipse (2cm and 1cm);
                \draw [line width=1.2pt, rotate=0] (-0.5,0) ellipse (1.5cm and 1.5cm);
            \end{tikzpicture}
        }
    } 
    & & \multirow{4}{*}{$\t$} \\
    & & one $A_3$, two $A_1$ & \\ 
    & & (one tacnode, two nodes)& \\ 
    & & & \\ 
    \hline 
    \multirow{4}{*}{ 
        \makecell[c]{ \\
            \begin{tikzpicture}
                \draw [style={double,double distance=2pt}] (0,0) -- (2,0);
                \draw [fill] (0,0) circle (3pt);
                \node at (0,-0.5) {$Q_1$};
                \draw [fill] (2,0) circle (3pt);
                \node at (2,-0.5) {$Q_2$};
            \end{tikzpicture}
        }
    }
    & 
    \multirow{4}{*}{
        \makecell[c]{
            \begin{tikzpicture}[scale=0.7]
                \draw [line width=1.2pt] (0,0) ellipse (2cm and 1.25cm);
                \draw [line width=1.2pt] (0,0) ellipse (1cm and 1.25cm);
            \end{tikzpicture}
        }
    }
    & & \multirow{4}{*}{$\tt$} \\
    & & two $A_3$ & \\ 
    & & (two tacnodes)& \\ 
    & & & \\ 
    \hline
    \multirow{4}{*}{ 
        \makecell[c]{ \\
            \begin{tikzpicture}
                \draw [style={decorate, decoration=snake}] (0,0) -- (2,0);
                \draw (0,0) -- (2,0);
                \draw [fill] (0,0) circle (3pt);
                \node at (0,-0.5) {$Q_1$};
                \draw [fill] (2,0) circle (3pt);
                \node at (2,-0.5) {$Q_2$};
            \end{tikzpicture}
        }
    }
    & 
    \multirow{4}{*}{
        \makecell[c]{
            \begin{tikzpicture}[scale=0.7]
                    \draw [line width=1.2pt, rotate around={10:(0.5,0)}] (2,0) ellipse (2cm and 1cm);
                    \draw [line width=1.2pt, rotate around={-10:(0.5,0)}] (2,0) ellipse (2cm and 1cm);
            \end{tikzpicture}
        }
    }
    & 
    \multirow{4}{*}{one $A_5$, one $A_1$}
    & 
    \multirow{4}{*}{$\a{5}$} \\
    & & & \\ 
    & & & \\ 
    & & & \\ 
    \hline
    \multirow{4}{*}{ 
        \makecell[c]{   \\
            \begin{tikzpicture}
                \draw [line width = 2.5pt] (0,0) -- (2,0);
                \draw [fill] (0,0) circle (3pt);
                \node at (0,-0.5) {$Q_1$};
                \draw [fill] (2,0) circle (3pt);
                \node at (2,-0.5) {$Q_2$};
            \end{tikzpicture}
        }
    }
    & 
    \multirow{4}{*}{
        \makecell[c]{
            \begin{tikzpicture}[scale=0.7]
                \draw [line width=1.2pt] (0,0) ellipse (2cm and 1.25cm);
                \draw [line width=1.2pt] (-0.5,0) ellipse (1.5cm and 1cm);
            \end{tikzpicture}
        }
    } 
    & 
    \multirow{4}{*}{one $A_7$}
    &
    \multirow{4}{*}{$\a{7}$} \\
    & & & \\ 
    & & & \\ 
    & & & \\ \hline
\end{tabular}

\caption{Notation for pairs of smooth conics.}
\label{tab: notation}

\end{table}

\begin{example}
    In this arrangement of three conics, given by
    \begin{align*}
        & Q_1: X^2 + Y^2 - Z^2 = 0, \\
        & Q_2: 4X^2 + Y^2 - Z^2 = 0, \\
        & Q_3: 4(X+Z)^2 - 6(XZ+Z^2) + 3Y^2 = 0,
    \end{align*}
    we have two $A_1$ singularities, three $A_3$ singularities and one $A_7$ singularity. This configuration splits into three pairs as follows:

    \begin{center}
    \begin{tikzpicture}[scale=1]
    \draw [line width=1.2pt] (0,-0.5) ellipse (1.2cm and 1.2cm);
    \draw [line width=1.2pt] (0,-0.5) ellipse (0.6cm and 1.2cm);
    \draw [line width=1.2pt] (-0.3,-0.5) ellipse (0.9cm and 1cm);
    
    \draw [line width=1.2pt] (0,-4) ellipse (1.2cm and 1.2cm);
    \draw [line width=1.2pt] (0,-4) ellipse (0.6cm and 1.2cm);
    \draw [dashed] (-0.3,-4) ellipse (0.9cm and 1cm);
    
    \draw [line width=1.2pt] (-3,-4) ellipse (1.2cm and 1.2cm);
    \draw [dashed] (-3,-4) ellipse (0.6cm and 1.2cm);
    \draw [line width=1.2pt] (-3.3,-4) ellipse (0.9cm and 1cm);
    
    \draw [dashed] (3,-4) ellipse (1.2cm and 1.2cm);
    \draw [line width=1.2pt] (3,-4) ellipse (0.6cm and 1.2cm);
    \draw [line width=1.2pt] (2.7,-4) ellipse (0.9cm and 1cm);
    
    \draw [-Stealth] (0,-2) -- (0,-2.6);
    \draw [-Stealth] (-1,-1.9) -- (-1.7,-2.6);
    \draw [-Stealth] (1,-1.9) -- (1.7,-2.6);
    \end{tikzpicture}
    \end{center}

    Using the notation from Table \ref{tab: notation}, we write this decomposition as $\a{7} + \tt + \t$.
\end{example}

\section{Free configurations of three conics with only ADE singularities}

Here is the list of all possible ADE singularities in configurations of three smooth plane conics:

\begin{enumerate}
    \item $A_1$ -- node, in which two conics intersect transversally,
    \vspace{-5pt}
    \item $A_3$ -- tacnode, in which two conics are tangent with multiplicity $2$,
    \vspace{-5pt}
    \item $A_5$ -- point in which two conics are tangent with multiplicity $3$,
    \vspace{-5pt}
    \item $A_7$ -- point in which two conics are tangent with multiplicity $4$,
    \vspace{-5pt}
    \item $D_4$ -- point where three conics intersect pairwise transversally,
    \vspace{-5pt}
    \item $D_6$ -- two conics $Q_1$ and $Q_2$ have an $A_3$ singularity and the third conic passes transversally through that singularity,
    \vspace{-5pt}
    \item $D_8$ -- two conics $Q_1$ and $Q_2$ have an $A_5$ singularity and the third conic passes transversally through that singularity,
    \vspace{-5pt}
    \item $D_{10}$ -- two conics $Q_1$ and $Q_2$ have an $A_7$ singularity and the third conic passes transversally through that singularity.
\end{enumerate}

Table \ref{table: ADES} shows all the notations and values of the Milnor number $\mu_p$ (which in this case is equal to the Tjurina number) and the multiplicities of all the singularities on this list. The values and the notation from this table will be used in the following sections of the paper.

\begin{table}
    \centering
    \begin{tabular}{c|c|c|c}
     & Number of & & \\
     Singularity & occurrences & $\mu_{p}$ & Multiplicity \\ \hline
    $A_{1}$ & $n_2$ & 1 & 1 \\
    $A_{3}$ & $t_3$ & 3 & 2 \\
    $D_{4}$ & $n_3$ & 4 & 3 \\
    $A_{5}$ & $t_5$ & 5 & 3 \\
    $D_{6}$ & $d_6$ & 6 & 4 \\
    $A_{7}$ & $t_7$ & 7 & 4 \\
    $D_{8}$ & $d_8$ & 8 & 5 \\
    $D_{10}$ & $d_{10}$ & 10 & 6
    \end{tabular}
    \caption {ADE singularities.}
    \label{table: ADES}
\end{table}

Let us recall some definitions first, based on \cite{Pokora2023}. The germ $(C, p)$ is weighted homogeneous of type $(w_1, w_2; 1)$ with $0 < w_j \leq \frac{1}{2}$ if there are local analytic coordinates $y_1, y_2$ centered at $p = (0,0)$ and a polynomial $g(y_1, y_2) = \sum\limits_{u,v}c_{u,v}y_1^uy_2^v$ with $c_{u,v} \in \mathbb{C}$, where the sum is over all pairs $(u,v) \in \mathbb{N}^2$ with $uw_1 + vw_2 = 1$. Using this description, the Arnold exponent (a log canonical threshold) of $p$ can be defined as
$$\alpha_p = w_1 + w_2.$$
If $\mathcal{C} = \{Q_1, Q_2, Q_3\}$ is a configuration of three smooth conics with only ${\rm ADE}$ singularities, then by  \cite[Theorem 2.1]{DIMCASERNESI} we know that  we have $d_1\geq \alpha_{\mathcal{C}}\cdot 6-2$ where  $\alpha_{\mathcal{C}}$  is the minimum
 of the Arnold exponents $\alpha_p$ of the singular points $p$ of $\mathcal{C}$. In the case of a configuration with only ${\rm ADE}$ singularities, we have $\alpha_{\mathcal{C}} > \frac{1}{2}$ (see \cite[Corollary 7.45 and inequality 7.46]{dimca2013topics}), so we have $d_1 > \frac{1}{2} \cdot 6 - 2 = 1$. Since $d_1 < \frac{d}{2}$ for free curves, we have $d_1 = 2$. Therefore by Theorem \ref{duPles} we get $\tau(\mathcal{C}) = 19$. Hence, using the values of $\mu_p$ and notation from the table above, we obtain the following system of equations:
\begin{equation}\label{system_equations}
\begin{cases}
    n_{2} + 3t_{3} + 4n_{3} + 5t_{5} + 6d_{6} + 7t_{7} + 8d_{8} + 10d_{10} = \tau(\mathcal{C}) = 19,\\
    n_{2} + 2t_{3} + 3n_{3} + 3t_{5} + 4d_{6} +  4t_{7} + 5d_{8} + 6d_{10} = 4\cdot\binom{3}{2} = 12,
\end{cases}
\end{equation}
where the second equation refers to the well-known naive combinatorial counting of the intersections among conics. By finding all solutions of \eqref{system_equations} in nonnegative integers, we obtain the following result:

\begin{thm}\label{thm: weakcombADE} A configuration of three smooth conics with only ${\rm ADE}$ singularities is free if and only if it has one of the following weak combinatorics:
 {\small
\begin{align*}
(n_{2}, t_{3}, n_{3}, t_{5}, d_{6}, t_{7}, d_{8}, d_{10}) \in \{(0, 3, 0, 2, 0, 0, 0, 0), (1, 1, 0, 3, 0, 0, 0, 0), \\
(0, 0, 1, 3, 0, 0, 0, 0), (0, 1, 0, 2, 1, 0, 0, 0), (0, 4, 0, 0, 0, 1, 0, 0), (1, 2, 0, 1, 0, 1, 0, 0), \\
(0, 1, 1, 1, 0, 1, 0, 0), (2, 0, 0, 2, 0, 1, 0, 0), (0, 2, 0, 0, 1, 1, 0, 0), (1, 0, 0, 1, 1, 1, 0, 0), \\
(0, 0, 0, 0, 2, 1, 0, 0), (2, 1, 0, 0, 0, 2, 0, 0), (1, 0, 1, 0, 0, 2, 0, 0), (0, 2, 0, 1, 0, 0, 1, 0), \\
(1, 0, 0, 2, 0, 0, 1, 0), (0, 0, 0, 1, 1, 0, 1, 0), (1, 1, 0, 0, 0, 1, 1, 0), (0, 0, 1, 0, 0, 1, 1, 0), \\
(0, 1, 0, 0, 0, 0, 2, 0), (0, 3, 0, 0, 0, 0, 0, 1), (1, 1, 0, 1, 0, 0, 0, 1), (0, 0, 1, 1, 0, 0, 0, 1), \\
(0, 1, 0, 0, 1, 0, 0, 1), (2, 0, 0, 0, 0, 1, 0, 1), (1, 0, 0, 0, 0, 0, 1, 1)\}.
\end{align*}
}
\end{thm}

Now we introduce some simple lemmas, which we will use later to show that most of these weak combinatorics cannot be realized geometrically over the complex numbers.

\begin{lemma}\label{split_a7}
    If a configuration of smooth conics with only ADE singularities has $t_7$ singularities of type $A_7$ and $d_{10}$ singularities of type $D_{10}$, then by splitting the configuration into pairs of conics, we obtain $t_7 + d_{10}$ pairs of type $\a{7}$.
\end{lemma}

\begin{proof}
    This is a straightforward observation resulting from the analysis of the list of singularities in arrangements of smooth conics at the beginning of this section. Namely, a pair of type $\a{7}$ can only be obtained from an $A_7$ singularity or from a $D_{10}$ singularity.
\end{proof}

\begin{lemma}\label{split_a5}
    If a configuration of smooth conics with only ADE singularities has $t_5$ singularities of type $A_5$ and $d_8$ singularities of type $D_8$, then by splitting the configuration into pairs of conics, we obtain $t_5 + d_8$ pairs of type $\a{5}$.    
\end{lemma}

\begin{proof}
    Analogous to the proof of Lemma \ref{split_a7}.
\end{proof}


\begin{lemma}\label{split_a3}
    If a configuration of smooth conics with only ADE singularities has $t_3$ singularities of type $A_3$ and $d_6$ singularities of type $D_6$, then by splitting the configuration into pairs of conics, we obtain $\ell_1$ pairs of type $\t$ and $\ell_2$ pairs of type $\tt$, with $\ell_1 + 2\ell_2 = t_3 + d_6$.
\end{lemma}

\begin{proof}
    This equality is a consequence of counting the $A_3$ singularities.
\end{proof}

\begin{lemma} \label{split_nod}
    If a configuration of three smooth conics with only ADE singularities has at least one singularity of type $A_7$ or $D_{10}$, then it does not contain singularities of type $D_4$, $D_6$ or $D_8$.
\end{lemma}

\begin{proof}
    If a configuration has a singularity of type $A_7$ or $D_{10}$, then there exists a pair of conics which intersect only at one point (with multiplicity $4$). However, a singularity of type $D_4$, $D_6$ or $D_8$ requires three conics to intersect and each intersection index is smaller than $4$, which is a contradiction.
\end{proof}

\begin{lemma} \label{split_nod2}
    If a configuration of three smooth conics with only ADE singularities has more than one pair without nodes (i.e.\ $\a{7}$ or $\tt$), then it does not contain singularities of type $D_6$, $D_8$ or $D_{10}$. 
\end{lemma}

\begin{proof}
    In order for the singularity of type $D_6$, $D_8$ or $D_{10}$ to exist, there must be two different pairs of conics with $A_1$ singularities. However, if a configuration has more than one pair without nodes, there exists at most one pair with $A_1$ singularities, a contradiction.
\end{proof}
\begin{lemma} \label{split_d8}
    If a configuration of three smooth conics with only ADE singularities has at least one singularity of type $D_8$ then each pair of conics has to contain nodes (i.e.\ $\n$, $\t$ or $\a{5}$).
\end{lemma}

\begin{proof}
    The $D_8$ singularity splits into three pairs, where two pairs have nodes and one pair has the $A_5$ singularity (i.e. $\a{5}$). A pair of type $\a{5}$ also contains a node.
\end{proof}

\begin{proposition}\label{prop: first}
    The following weak combinatorics:
    $$(n_{2}, t_{3}, n_{3}, t_{5}, d_{6}, t_{7}, d_{8}, d_{10}) \in \{(0,3,0,2,0,0,0,0), (1,1,0,3,0,0,0,0)\}$$
    cannot be realized over $\mathbb{C}$ by three smooth conics.
\end{proposition}

\begin{proof}
    From Lemmas \ref{split_a5} and, \ref{split_a3} we obtain at least $4$ pairs in the decomposition in each case, which is a contradiction.
\end{proof}

\begin{proposition}
    The following weak combinatorics
\begin{align*}
(n_{2}, t_{3}, n_{3}, t_{5}, d_{6}, t_{7}, d_{8}, d_{10}) \in \{(0, 1, 1, 1, 0, 1, 0, 0), (0, 2, 0, 0, 1, 1, 0, 0), \\
(1, 0, 0, 1, 1, 1, 0, 0), (0, 0, 0, 0, 2, 1, 0, 0), (1, 0, 1, 0, 0, 2, 0, 0), (1, 1, 0, 0, 0, 1, 1, 0), \\(0, 0, 1, 0, 0, 1, 1, 0), (0, 0, 1, 1, 0, 0, 0, 1), (0, 1, 0, 0, 1, 0, 0, 1), (1, 0, 0, 0, 0, 0, 1, 1))\}.
\end{align*}
    cannot be realized over $\mathbb{C}$ by smooth conics.
\end{proposition}

\begin{proof}
    Follows from Lemma \ref{split_nod}.  
\end{proof}

\begin{proposition}
    The following weak combinatorics
$$(n_{2}, t_{3}, n_{3}, t_{5}, d_{6}, t_{7}, d_{8}, d_{10}) \in \{(2, 0, 0, 0, 0, 1, 0, 1), (0, 3, 0, 0, 0, 0, 0, 1)\}$$
    cannot be realized over $\mathbb{C}$ by smooth conics.
\end{proposition}

\begin{proof}
    Follows from Lemma \ref{split_nod2}.    
\end{proof}

\begin{proposition}
    The configuration of three smooth conics with the weak combinatoric
    $$(n_{2}, t_{3}, n_{3}, t_{5}, d_{6}, t_{7}, d_{8}, d_{10}) = (0, 2, 0, 1, 0, 0, 1, 0)$$
    cannot be realized over $\mathbb{C}$.
\end{proposition}

\begin{proof}
    Suppose that such a configuration exists. Then it splits into $2$ pairs of type $\a{5}$ and one pair of type $\tt$. Since $\tt$ does not contain a node, this weak combinatoric cannot be realized from Lemma \ref{split_d8}.
\end{proof}

\begin{proposition}
    The configuration of three smooth conics with the weak combinatoric
$$(n_{2}, t_{3}, n_{3}, t_{5}, d_{6}, t_{7}, d_{8}, d_{10}) = (0, 0, 0, 1, 1, 0, 1, 0)$$
    cannot be realized over $\mathbb{C}$.
\end{proposition}

\begin{proof}
    Suppose that such a configuration exists. Then it splits into $2$ pairs of type $\a{5}$ and one pair of type $\t$. The singularity $D_6$ is formed by joining the singularity $A_3$ with two nodes. These nodes must come from two $\a{5}$ pairs. The singularity $D_8$ is formed by joining the singularity $A_5$ with two nodes. Since we already used the nodes from both $\a{5}$ pairs, the only remaining ones are from the $\t$ pair, but these two nodes are clearly distinct, a contradiction.
\end{proof}

Before we proceed, we introduce the following two facts which will be used in our proofs.

\begin{fact}\label{fact:all} Let $Q_1$ and $Q_2$ be two distinct smooth conics.
\begin{enumerate}[label=\alph*)]
    \item \label{fact:all_a5} If $Q_1$ and $Q_2$ form the $\a{5}$ arrangement (i.e.\ with singularities $A_5$ and $A_1$), then the equation of $Q_2$ can be written in the form
    $$Q_2: \lambda Q_1 + \ell_1\ell_2 = 0,$$
    where $\lambda \in \mathbb{C} \setminus \{0\}$, $\ell_1$ is the equation of a line tangent to $Q_1$ at the $A_5$ singularity and $\ell_2$ is the equation of a line passing through both singularities.
    \item \label{fact:all_t} If $Q_1$ and $Q_2$ form the $\t$ arrangement (i.e.\ with one $A_3$ singularity and two nodes), then the equation of $Q_2$ can be written in the form
    $$Q_2: \lambda Q_1 + \ell_1\ell_2 = 0,$$
    where $\lambda \in \mathbb{C} \setminus \{0\}$, $\ell_1$ is the equation of a line tangent to $Q_1$ at the $A_3$ singularity and $\ell_2$ is the equation of a line passing through both nodes.
    \item \label{fact:all_t2} If $Q_1$ and $Q_2$ form the $\tt$ arrangement (i.e.\ with two $A_3$ singularities), then the equation of $Q_2$ can be written in the form
    $$Q_2: \lambda Q_1 + \ell_1\ell_2 = 0,$$
    where $\lambda \in \mathbb{C} \setminus \{0\}$ and $\ell_1$, $\ell_2$ are the equation of lines tangent to $Q_1$ at the respective singularities.
    \item \label{fact:all_a7} If $Q_1$ and $Q_2$ form the $\a{7}$ arrangement (i.e.\ with one $A_7$ singularity), then the equation of $Q_2$ can be written in the form
    $$Q_2: \lambda Q_1 + \ell^2 = 0,$$
    where $\lambda \in \mathbb{C} \setminus \{0\}$ and $\ell$ is the equation of a line tangent to $Q_1$ at the $A_7$ singularity.
\end{enumerate}
\end{fact}

\begin{proposition}
    The following weak combinatoric
$$(n_{2}, t_{3}, n_{3}, t_{5}, d_{6}, t_{7}, d_{8}, d_{10}) = (0, 4, 0, 0, 0, 1, 0, 0)$$
    \noindent cannot be realized over $\mathbb{C}$ by smooth conics.
\end{proposition}

\begin{proof}
Suppose that there exists a configuration admitting this weak combinatoric. The only possible decomposition into pairs is $\a{7} + 2 \cdot \tt$. Using the method presented in \cite{PHD} (see Table \ref{tab: notation}), this configuration of conics can be represented by the following graph.

\begin{center}
\begin{tikzpicture}
\draw [style={double,double distance=2pt}] (0,0) -- (2,0);
\draw [line width = 2.5pt] (0,0) -- (1,1.5);
\draw [style={double,double distance=2pt}] (1,1.5) -- (2,0);
\draw [fill] (0,0) circle (3pt);
\node at (-0.5,0) {$Q_1$};
\draw [fill] (2,0) circle (3pt);
\node at (2.5,0) {$Q_2$};
\draw [fill] (1,1.5) circle (3pt);
\node at (1,2) {$Q_3$};
\end{tikzpicture}
\end{center}
By \cite[Proposition 3]{Megyesi2000}, we may assume that
$$Q_1: X^2 + Y^2 - Z^2 = 0, \quad Q_2: \frac{1}{q^2}X^2 + Y^2 - Z^2 = 0,$$
for some $q \in \mathbb{C}$, $q \notin \{0, 1, -1\}$. By \cite[Remark 4.3.2]{PHD} we know that $Q_1$, $Q_2$ have following parametrizations:
\begin{align*}
Q_1 &= \{[2st : t^2 - s^2 : t^2 + s^2 ] \mid [s:t] \in \mathbb{P}^1 \} \\
Q_2 &= \{[2qst  : t^2-s^2:t^2 + s^2] \mid  [s:t] \in \mathbb{P}^1 \}
\end{align*}
Let us assume that $Q_3$ is given by $a_1X^2 + a_2Y^2 + a_3Z^2 + a_4XY + a_5YZ + a_6XZ = 0$. Since $Q_1$ and $Q_2$ intersect in $[0:1:\pm 1]$ and $Q_3$ cannot contain these points, then by parametrizations of $Q_1$ and $Q_2$ we obtain that the intersection points of $Q_3$ with $Q_1$ and $Q_2$ must be of the form
$$[2u : u^2 - 1 : u^2 + 1] \text{ and } [2qv : v^2 - 1 : v^2 + 1],$$
respectively, where $u \neq 0$ and $v \neq 0$. By substituting these points into the equation of $Q_3$, we obtain
\begin{align*}
f_{13}(u) &= (a_2 + a_3 + a_5)u^4 + 2(a_4 + a_6)u^3 + (4a_1 - 2a_2 + 2a_3)u^2 \\
&+ 2(-a_4 + a_6)u + (a_2 + a_3 - a_5) = 0
\end{align*}
and
\begin{align*}
f_{23}(v) &= (a_2 + a_3 + a_5)v^4 + 2q(a_4 + a_6)v^3 + (4a_1q^2 - 2a_2 + 2a_3)v^2 \\
&+ 2q(-a_4 + a_6)v + (a_2 + a_3 - a_5) = 0
\end{align*}
By the intersection behaviour of $Q_3$ with $Q_1$ and $Q_2$, we have $f_{13}(u) = A(u - \lambda)^4$ and $f_{23}(v) = B(v - \mu_1)^2(v - \mu_2)^2$ for some $A, B,\lambda, \mu_1, \mu_2 \in \mathbb{C}$, with $A \neq 0$, $B \neq 0$ and $\mu_1 \neq \mu_2$. By comparing the coefficients, we obtain the following system of equations:
\begin{center}
    $\left\{
\begin{array}{lr}
a_2 + a_3 + a_5  = A & \text{(I)} \\
2(a_4 + a_6)  = -4A\lambda & \text{(II)} \\
4a_1 - 2a_2 + 2a_3 = 6A\lambda^2 & \text{(III)} \\
2(-a_4 + a_6) = -4A\lambda^3 & \text{(IV)} \\
a_2 + a_3 - a_5 = A\lambda^4 & \text{(V)} \\
a_2 + a_3 + a_5 = B & \text{(VI)} \\
2q(a_4 + a_6)  = -2B(\mu_1 + \mu_2) & \text{(VII)} \\
4a_1q^2 - 2a_2 + 2a_3  = B(\mu_1^2 + 4\mu_1\mu_2 + \mu_2^2) & \text{(VIII)} \\
2q(-a_4 + a_6)  = -2B\mu_1\mu_2(\mu_1 + \mu_2) & \text{(IX)} \\
a_2 + a_3 - a_5  = B\mu_1^2\mu_2^2 & \text{(X)} \\
\end{array}
\right.$
\end{center}
If $\lambda = 0$, then from (V) we get $a_2 + a_3 - a_5 = 0$, but this implies $[0:1:-1] \in Q_3$ which is a contradiction, hence we know that $\lambda \neq 0$. 

From (I) and (VI) we obtain $A=B$. Then (II) and (VII) implies $-4Aq\lambda = -2B(\mu_1 + \mu_2)$, which is equivalent to 
\begin{align*}
    &\mu_1 + \mu_2 = 2\lambda q. \tag{*}
\end{align*}
Now, from (IV), (IX) and (*), we obtain
\begin{align*}
-&4Aq\lambda^3 = -2B \mu_1 \mu_2 (\mu_1 + \mu_2) \\ 
-&4q\lambda^3 = -4\mu_1 \mu_2 \lambda q \\ 
&\lambda^2 = \mu_1 \mu_2. \tag{**}
\end{align*}
By adding (I) and (V), we get
\begin{align*}
    a_2 + a_3 &=  \frac{1}{2}A(1 + \lambda ^4). \tag{***}
\end{align*}
Now, after subtracting (I) and (V)  we obtain $a_5 = \frac{1}{2}A(1-\lambda^4)$. Similarly, from (II) and (IV) we get
\begin{align*}
    & a_6 =  -A\lambda (1 + \lambda^2), \\ 
    & a_4 =  -A\lambda (1 - \lambda^2),
\end{align*}
respectively. From (III) we get $-2a_2 + 2a_3 = 6A \lambda^2 - 4a_1$. Now, by substituting it into (VIII) and using (*) and (**), we obtain
\begin{align*}
    &4a_1q^2 - 2a_2 +2a_3 = A((\mu_1 + \mu_2)^2 + 2\mu_1 \mu_2) \\ 
    & a_1(q^2 - 1) = A\lambda^2(q^2-1) \\
    & a_1 = A\lambda^2.
\end{align*}
Substituting $a_1 = A\lambda^2$ in (III), we get
\begin{align*}
   -a_2 + a_3 = A\lambda^2. \tag{****}
\end{align*}
Now, by (***) and (****) we can find formulas for $a_2$ and $a_3$:
\begin{align*}
    &a_2 = \frac{1}{4}A(1-\lambda^2)^2 \\
    &a_3 =  \frac{1}{4}A(1+\lambda^2)^2 .
\end{align*}
Having determined $a_1,\ldots,a_6$ in terms of $A$ and $\lambda$ we can express $Q_3$ by  
\begin{align*}
    &Q_3: A\lambda^2 \cdot X^2 + \frac{1}{4}A(1-\lambda^2)^2 \cdot Y^2 + \frac{1}{4}A(1+\lambda^2)^2 \cdot Z^2 
    - A\lambda(1-\lambda^2) \cdot XY  \\ & \quad + \frac{1}{2}A(1-\lambda^4) \cdot YZ - A\lambda(1 + \lambda^2) XZ \\ 
    & \quad = A\left[ -\lambda \cdot X + \frac{1}{2}(1-\lambda^2) \cdot Y + \frac{1}{2}(1 + \lambda^2) \cdot Z\right]^2 = 0,
\end{align*} which is a contradiction with irreducibility of $Q_3$.  
\end{proof}

\begin{proposition}
    The following weak combinatoric
    $$(n_{2}, t_{3}, n_{3}, t_{5}, d_{6}, t_{7}, d_{8}, d_{10}) = (1, 2, 0, 1, 0, 1, 0, 0)$$
    \noindent cannot be realized over $\mathbb{C}$ by smooth conics.
\end{proposition}

\begin{proof}
    Suppose that such a configuration of three conics $Q_1$, $Q_2$ and $Q_3$ exists and it is represented by the following graph:
    
    \begin{center}
    \begin{tikzpicture}
    \draw [line width = 2.5pt] (0,0) -- (2,0);
    \draw [style={decorate, decoration=snake}] (0,0) -- (1,1.5);
    \draw (0,0) -- (1,1.5);
    \draw [style={double,double distance=2pt}] (1,1.5) -- (2,0);
    \draw [fill] (0,0) circle (3pt);
    \node at (-0.5,0) {$Q_1$};
    \draw [fill] (2,0) circle (3pt);
    \node at (2.5,0) {$Q_2$};
    \draw [fill] (1,1.5) circle (3pt);
    \node at (1,2) {$Q_3$};
    \end{tikzpicture}
    \end{center}
    By \cite[Proposition 4.2.1]{PHD}, we may assume that
    $$Q_1: X^2 - YZ = 0, \quad Q_2: X^2 + aZ^2 - YZ = 0,$$
    for some $a \in \mathbb{C}\setminus\{0\}$. These two conics meet at the point $[0 : 1 : 0]$ with multiplicity $4$ and they have the following parametrizations:
    \begin{align*}
        Q_1 &= \{[uv : v^2 : u^2] \mid [u : v] \in \mathbb{P}^1\}, \\
        Q_2 &= \{[st : t^2 + as^2: s^2] \mid  [s : t] \in \mathbb{P}^1\}.
    \end{align*}
    Suppose such a conic $Q_3$ exists. Since $[0 : 1 : 0] \notin Q_1 \cap Q_3$, then $Q_3$ intersects $Q_1$ at the points $[p : p^2 : 1]$ and $[q : q^2 : 1]$ (where $p \neq q$) with multiplicities $3$ and $1$, respectively. Similarly, $Q_3$ intersects $Q_1$ at the points $[m : m^2 + a : 1]$ and $[n : n^2 + a : 1]$ (where $m \neq n$), with multiplicity $2$. In addition, the line $\ell_1: 2pX - Y - p^2Z = 0$ is tangent to $Q_1$ at $[p : p^2 : 1]$ and the line $\ell_2: (p+q)X - Y - pqZ = 0$ passes through the intersection points $[p : p^2 : 1]$ and $[q : q^2 : 1]$. Therefore by Fact \ref{fact:all}~\ref{fact:all_a5} the equation of $Q_3$ must be of the form
    \begin{align*}
        &Q_3: \lambda Q_1 + \ell_1\ell_2 = (\lambda + 2p(p+q))X^2 + Y^2 + p^3qZ^2 - (3p+q)XY \\
        & \quad + (p(p+q) - \lambda)YZ - p^2(p + 3q)XZ = 0,
    \end{align*}
    for some $\lambda \in \mathbb{C}\setminus\{0\}$. Substituting the affine parametrization $x = t$, $y = t^2 + a$, $z = 1$ of $Q_2$ into the above equation, we obtain
    \begin{align*}
        f_{23}(t) &= t^4 - (3p+q)t^3 + (3p(p+q) + 2a)t^2 - (p^2(p + 3q) + a(3p + q))t \\
        & + (p^3q + a^2 + ap(p + q) - a\lambda) = 0.
    \end{align*}
    On the other hand, by the intersection behavior of $Q_2$ and $Q_3$, $f_{23}(t)$ must be of the form
    \begin{align*}
        f_{23}(t) = (t - m)^2(t - n)^2 &= t^4 - 2(m+n)t^3 + (m^2 + 4mn + n^2)t^2 \\
        & - 2mn(m+n)t + m^2n^2 = 0.
    \end{align*}
    Comparing these two equations, we obtain    $$\left\{\begin{array}{l}
        2(m + n) = 3p + q \\
        m^2 + 4mn + n^2 = 3p(p + q) + 2a \\
        2mn(m + n) = p^2(p + 3q) + a(3p + q) \\
        m^2n^2 = p^3q + a^2 + ap(p + q) - a\lambda
    \end{array}\right.$$ 
    From the first $2$ equations, we obtain
    $$m+n = \frac{3}{2}p + \frac{1}{2}q, \quad mn = \frac{3}{8}p^2+\frac{3}{4}pq - \frac{1}{8}q^2 + a.$$
    Then, by substituting these equalities into the third equation and simplifying, we finally obtain that $(p-q)^3 = 0$, which contradicts our earlier assumption that $p \neq q$.
\end{proof}


\begin{proposition}
    The following weak combinatoric
    $$(n_{2}, t_{3}, n_{3}, t_{5}, d_{6}, t_{7}, d_{8}, d_{10}) = (0, 1, 0, 2, 1, 0, 0, 0)$$
    \noindent cannot be realized over $\mathbb{C}$ by smooth conics.
\end{proposition}

\begin{proof}
    Suppose that such a configuration of three conics $Q_1$, $Q_2$ and $Q_3$ exists. It can be split into pairs $2 \cdot \a{5} + \tt$ and it can be represented by the following graph:
    \begin{center}
    \begin{tikzpicture}
    \draw [style={decorate, decoration=snake}] (0,0) -- (2,0);
    \draw (0,0) -- (2,0);
    \draw [style={decorate, decoration=snake}] (0,0) -- (1,1.5);
    \draw (0,0) -- (1,1.5);
    \draw [style={double,double distance=2pt}] (1,1.5) -- (2,0);
    \draw [fill] (0,0) circle (3pt);
    \node at (-0.5,0) {$Q_1$};
    \draw [fill] (2,0) circle (3pt);
    \node at (2.5,0) {$Q_2$};
    \draw [fill] (1,1.5) circle (3pt);
    \node at (1,2) {$Q_3$};
    \end{tikzpicture}
    \end{center}
    Using projective transformations, we can assume that the $D_6$ singularity is in $[0:0:1]$ and two $A_5$ singularities are in $[1:0:0]$ and $[0:1:0]$. The conic $Q_1$ passes through both $A_5$ singularities. Since it also has to pass through the $D_6$ singularity, it is necessarily given by
    $$Q_1: \; XY + aXZ + bYZ = 0$$
    for some $a, b \in \mathbb{C} \setminus \{0\}$. The lines tangent to $Q_1$ at $[1:0:0]$ and $[0:1:0]$ are $Y + aZ = 0$ and $X + bZ = 0$, respectively. The line passing through both intersection points of $Q_1$ and $Q_2$ is $Y = 0$, while the line passing through both intersection points of $Q_1$ and $Q_3$ is $X = 0$. Therefore, by Fact \ref{fact:all}~\ref{fact:all_a5}, the equations of the remaining conics are:
    \begin{align*}
        Q_2: \;\lambda_2Q_1 + (Y+aZ)\cdot Y = Y^2 + \lambda_2XY + \lambda_2aXZ + (\lambda_2b + a)YZ = 0, \\
        Q_3: \;\lambda_3Q_1 + (X+bZ)\cdot X = X^2 + \lambda_3XY + (\lambda_3a + b)XZ + \lambda_3bYZ= 0,
    \end{align*}
    for some $\lambda_2,\lambda_3 \in \mathbb{C} \setminus \{0\}.$ We want conics $Q_2$ and $Q_3$ to be tangent to each other at two points and one of them is $[0:0:1]$. The lines tangent to $Q_2$ and $Q_3$ at this point are
    $$\lambda_2aX + (\lambda_2b+a)Y = 0 \qquad \text{and} \qquad(\lambda_3a+b)X + \lambda_3bY = 0,$$
    respectively. For these lines to be identical, we need
    \begin{equation}\label{det_eq}
        \det\begin{bmatrix}
            \lambda_2a & \lambda_2b + a \\ \lambda_3a + b & \lambda_3b
        \end{bmatrix} = -\lambda_2b^2-\lambda_3a^2-ab = 0.        
    \end{equation}
Moreover, the coordinates of points of intersection of $Q_2$ and $Q_3$ satisfy
$$\left \{
\begin{array}{lr}
Y^2 + \lambda_2XY + \lambda_2aXZ + (\lambda_2b + a)YZ = 0, & \text{(I)} \\
X^2 + \lambda_3XY + (\lambda_3a + b)XZ + \lambda_3bYZ = 0, & \text{(II)}
\end{array} \right.$$
We can use this system of equations to determine the second point of tangency. Note that if $X = 0$, then $Y = 0$ and we obtain the point $[0 : 0 : 1]$, so we may assume $X = 1$. By multiplying equations (I) and (II) by $\lambda_3b$ and $(\lambda_2b + a)$, respectively, and subtracting them, we get
$$\lambda_3bY^2 + (\lambda_2\lambda_3b - \lambda_3(\lambda_2b + a))Y + (\lambda_2\lambda_3ab - (\lambda_3a + b)(\lambda_2b + a))Z - (\lambda_2b + a) = 0,$$
which by \eqref{det_eq} further reduces to
$$\lambda_3bY^2 - \lambda_3aY - (\lambda_2b + a) = 0.$$
Since we only have one more point of intersection of $Q_2$ and $Q_3$ besides $[0 : 0 : 1]$, the discriminant of the above quadratic equation must be zero. However, by \eqref{det_eq} this discriminant is equal to $-3\lambda_3^2a^2$, which is a contradiction since $\lambda_3 \neq 0$ and $a \neq 0$.
\end{proof}

\begin{proposition}\label{prop: base}
Any free configuration of three conics $Q_1$, $Q_2$, $Q_3$ with weak combinatorics
$$(n_{2}, t_{3}, n_{3}, t_{5}, d_{6}, t_{7}, d_{8}, d_{10}) = (0, 0, 1, 3, 0, 0, 0, 0)$$
is projectively equivalent to the conics
\begin{align*}
    &Q_1: X^2 - YZ = 0, \\
    &Q_2: 3X^2 + Y^2 + 3XY - 3YZ = 0, \\
    &Q_3: X^2 + 27Z^2 + 9XZ - YZ = 0.
\end{align*}
\end{proposition}

\begin{proof}
The decomposition into pairs is $3 \cdot \a{5}$. By \cite[Proposition 4.2.1]{PHD} we may assume that
\begin{align*}
    &Q_1: X^2 - YZ = 0, \\
    &Q_2: X^2 + bY^2 + cXY - YZ = 0 
\end{align*}
for some $b,c \in \mathbb{C}, c \neq 0$. These conics intersect at points $[0 : 0 : 1]$ and $[-bc : c^2 : b^2]$ with multiplicities $3$ and $1$, respectively (and the second point will become a $D_4$ singularity). These two conics have the following parametrizations:
\begin{align*}
    &Q_1 = \{[uv : v^2 : u^2] \mid [u : v] \in \mathbb{P}^1\}, \\
    &Q_2 = \{[st : s^2 : t^2 + bs^2 + cst] \mid [s : t] \in \mathbb{P}^1\}.
\end{align*}
We can therefore write the coordinates of the remaining two singularities of type $A_5$ as $[q : 1 : q^2]$ and $[r : 1 : r^2 + br + c]$ for some $q, r \in \mathbb{C}$ with $q \neq -\frac{b}{c}$ and $r \neq -\frac{b}{c}$.

By Fact \ref{fact:all}~\ref{fact:all_a5}, we can write the equation defining $Q_3$ as $\lambda Q_1 + \ell_1\ell_2$, where $\lambda \neq 0$, $\ell_1$ is the equation of a line tangent to $Q_1$ at point $[q : 1 : q^2]$, and $\ell_2$ is the equation of a line passing through points $[q : 1 : q^2]$ and $[-bc : c^2 : b^2]$. In this case, $\ell_1$ is defined by equation $2qX - q^2Y - Z = 0$ and $\ell_2$ is defined by $(b-cq)X - bqY + cZ = 0$. We therefore obtain
\begin{align*}
    Q_3 & : \lambda(X^2 - YZ) + (2qX - q^2Y - Z)((b-cq)X - bqY + cZ) \\
    &= (2bq - 2cq^2 + \lambda)X^2 + (-3bq^2 + cq^3)XY + bq^3Y^2 \\
    & + (-b + 3cq)XZ + (bq - cq^2 - \lambda)YZ - cZ^2 = 0.
\end{align*}
On the other hand, we can write the equation for $Q_3$ in the same way, but this time considering its intersection with $Q_2$, i.e.\ in the form $\mu Q_2 + \ell_3\ell_4$, where $\mu \neq 0$ and the lines $\ell_3$ and $\ell_4$ are defined similarly as before. In this case, we obtain
\begin{align*}
    &\ell_3 : (2r + c)X + (b - r^2)Y - Z = 0, \\
    &\ell_4 : (b - cr - c^2)X - b(c+r)Y + cZ = 0,
\end{align*}
and therefore
\begin{align*}
    Q_3 & : \mu(X^2 + bY^2 + cXY - YZ) + ((2r + c)X + (b - r^2)Y - Z)((b - cr - c^2)X - b(c+r)Y + cZ) \\
    &= (\mu + bc + 2br - c^3 - 3c^2r - 2cr^2)X^2 + (\mu c + b^2 - 2bc^2 - 4bcr - 3br^2 + c^2r^2 + cr^3)XY\\
    & + (\mu b - b^2c - b^2r + bcr^2 + br^3)Y^2 + (-b + 2c^2 + 3cr)XZ + (-\mu + 2bc + br - cr^2)YZ - cZ^2 = 0.
\end{align*}
By comparing the coefficients in both equations, we obtain the following system of equations:
$$\left \{
\begin{array}{lr}
2bq - 2cq^2 + \lambda = \mu + bc + 2br - c^3 - 3c^2r - 2cr^2 & \text{(I)} \\
-3bq^2 + cq^3 = \mu c + b^2 - 2bc^2 - 4bcr - 3br^2 + c^2r^2 + cr^3 & \text{(II)}\\ 
bq^3 = \mu b - b^2c - b^2r + bcr^2 + br^3  & \text{(III)} \\ 
-b+3cq = -b + 2c^2 + 3cr & \text{(IV)} \\
bq - cq^2 - \lambda = -\mu + 2bc + br - cr^2 & \text{(V)}
\end{array} \right.$$
From (IV) we obtain $q = \frac{2}{3}c + r$. By substituting that into the remaining equations, we obtain $\lambda = \mu$ (equations (I) and (V)) and $r = -\frac{b}{c} - \frac{1}{3}c$ (equations (II) and (III)), which leads to $q = -\frac{b}{c} + \frac{1}{3}c$. Finally, from equation (III) we get $\mu = \lambda = -\frac{1}{27}c^3$. By substituting these values into any of the above equations defining $Q_3$, we finally obtain the following equations:
\begin{align}\label{eq1_param}
\begin{split}
    &Q_1: X^2 - YZ = 0, \\
    &Q_2: X^2 + bY^2 + cXY - YZ = 0, \\
    &Q_3: c^2(7c^4 - 54bc^2 + 108b^2)X^2 - c(c^2 - 12b)(c^2 - 3b)^2XY - b(c^2 - 3b)^3Y^2 \\
    & \quad - 27c^3(c^2 - 4b)XZ + c^2(2c^4 - 27bc^2 + 54b^2)YZ + 27c^4Z^2 = 0,
\end{split}
\end{align}
with $b, c \in \mathbb{C}$ and $c \neq 0$. The singularities in the obtained arrangement have the following coordinates:
\begin{itemize}
    \item $[0 : 0 : 1]$ -- an $A_5$ singularity, intersection of $Q_1$ and $Q_2$,
    \vspace{-6pt}
    \item $\big[-bc+\frac{1}{3}c^3 : c^2 : (b - \frac{1}{3}c^2)^2\big]$ -- an $A_5$ singularity, intersection of $Q_1$ and $Q_3$,
    \vspace{-6pt}
    \item $\big[-bc - \frac{1}{3}c^3 : c^2 : (b+\frac{1}{3}c^2)^2 - \frac{1}{3}c^4\big]$ -- an $A_5$ singularity, intersection of $Q_2$ and $Q_3$,
    \vspace{-6pt}
    \item $[-bc : c^2 : b^2]$ -- a $D_4$ singularity.
\end{itemize}
It can be checked that no three of these points are collinear, we can therefore map them to the points $[1 : 0 : 0]$, $[0 : 1 : 0]$, $[0 : 0 : 1]$ and $[1 : 1 : 1]$, respectively. If we denote this map by $M$, then
$$M^{-1} = \begin{bmatrix} 
0 & -9bc+3c^3 & -9bc-3c^3 \\ 
0 & 9c^2 & 9c^2 \\ 
c^4 & (3b-c^2)^2 & 9b^2+6bc^2-2c^4 
\end{bmatrix},$$
and therefore the conics are mapped to
\begin{align*}
    &Q_1': 3Z^2 - XY - XZ - YZ = 0, \\
    &Q_2': 3Y^2 - XY - XZ - YZ = 0, \\
    &Q_3': 3X^2 - XY - XZ - YZ = 0.
\end{align*}
Hence all arrangements given by \eqref{eq1_param} are projectively equivalent. In particular they are equivalent to \eqref{eq1_param} for $b = \frac{1}{3}$ and $c = 1$, which are the equations from the statement of our proposition.
\end{proof}

\begin{proposition}\label{prop:twocases}
Any free configuration of three conics $Q_1$, $Q_2$, $Q_3$ with weak combinatorics
$$(n_{2}, t_{3}, n_{3}, t_{5}, d_{6}, t_{7}, d_{8}, d_{10}) = (2, 0, 0, 2, 0, 1, 0, 0)$$
is projectively equivalent to the conics
\begin{align*}
    &Q_1: X^2 - YZ = 0, \\
    &Q_2: X^2 + i\sqrt{3}Z^2 - YZ = 0, \\
    &Q_3: 8X^2 - 3(i\sqrt{3}+1)Y^2 + 8i\sqrt{3}XY - 8YZ = 0.
\end{align*}
\end{proposition}

\begin{proof}
The statement and proof of the above theorem was already published in \cite{PHD}, however, the formula presented there contains mistakes. The correct formulas are the following:
\begin{align}\label{eq2_param}
\begin{split}
    &Q_1: X^2 - YZ = 0, \\
    &Q_2: X^2 + aZ^2 - YZ = 0, \\
    &Q_3: \left(\frac{(m+p)(m-p)^3}{a} + \frac{a(m+p)}{m-p} - \frac{2}{3}(m^2 - 8mp + p^2)\right)X^2 + Y^2 \\
    & \quad + \left(\frac{2ap^3}{3(m-p)} + 2mp^3 - p^4\right)Z^2 - \left(\frac{2a}{3(m-p)} + 2m + 2p\right)XY \\
    & \quad - \left(\frac{2ap^2}{m-p} + 6mp^2 - 2p^3\right)XZ - \left(\frac{(m+p)(m-p)^3}{a} + a - \frac{2}{3}(m^2 + mp + p^2)\right)YZ = 0,
\end{split}
\end{align}
with $a, m, p \in \mathbb{C}$, $m \neq p$ and $a^2 = -3(m-p)^4$. This arrangement has the following singularities:
\begin{itemize}
    \item $[0 : 1 : 0]$ -- an $A_7$ singularity, intersection of $Q_1$ and $Q_2$,
    \vspace{-6pt}
    \item $[p : p^2 : 1]$ -- an $A_5$ singularity, intersection of $Q_1$ and $Q_3$,
    \vspace{-6pt}
    \item $[m : m^2 + a : 1]$ -- an $A_5$ singularity, intersection of $Q_2$ and $Q_3$,
    \vspace{-6pt}
    \item $[q : q^2 : 1]$ -- an $A_1$ singularity, intersection of $Q_1$ and $Q_3$,
    \vspace{-6pt}
    \item $[n : n^2+a : 1]$ -- an $A_1$ singularity, intersection of $Q_2$ and $Q_3$,
    \end{itemize}
    with $q = 2m - p + \frac{2a}{3(m-p)}$ and $n = 2p - m + \frac{2a}{3(m-p)}.$ The first four of these points can be mapped into $[1:0:0]$, $[0:1:0]$, $[0:0:1]$ and $[1:1:1]$, respectively. If we denote this map by $M$, then
$$M^{-1} = \begin{bmatrix}
0 & (-2u - 3)p & (2u + 6)m \\
8(m-p)^2 & (-2u - 3)p^2 & (2u + 6)(m^2 + a) \\
0 & -2u - 3 & 2u + 6
\end{bmatrix},$$
with $a = u(m-p)^2$, either $u = i\sqrt{3}$ or $u = -i\sqrt{3}$. The conics are mapped to
\begin{align*}
    &Q_1': (3-9u)Z^2 + 7XY + (2u - 10)XZ + 7uYZ = 0, \\
    &Q_2': (21u - 42)Y^2 + 56XY + (16u-80)XZ + (84-28u)YZ = 0, \\
    &Q_3': 12X^2 + (u - 9)XY + 8uXZ - (9u + 3)YZ = 0.
\end{align*}
For both values of $u$ we obtain projectively equivalent arrangements. The matrix of the projective mapping from the arrangement for $u=-i\sqrt{3}$ to the arrangement for $u = i\sqrt{3}$ is
$$N = \begin{bmatrix}
6i\sqrt{3}+2 & 0 & 0 \\
0 & 0 & -16 \\
0 & -7 & 0
\end{bmatrix}.$$
Since all the considered arrangements are projectively equivalent, we can choose $m = 1$, $p = 0$ and $a = i\sqrt{3}$ in \eqref{eq2_param}, obtaining the equations from the statement of our proposition.
\end{proof}

\begin{remark}
    Since the parameters $a, m$ and $p$ have to satisfy the equality $a^2 = -3(m-p)^4$ with $m \neq p$, this arrangement cannot be realized over $\mathbb{R}$.
\end{remark}

\begin{proposition}\label{prop: a3d8d8}
Any free configuration of three conics $Q_1$, $Q_2$, $Q_3$ with weak combinatorics
$$(n_{2}, t_{3}, n_{3}, t_{5}, d_{6}, t_{7}, d_{8}, d_{10}) = (0, 1, 0, 0, 0, 0, 2, 0)$$
is projectively equivalent to the conics
\begin{align*}
    &Q_1: X^2 - YZ = 0, \\
    &Q_2: X^2 - 2XY - YZ = 0, \\
    &Q_3: X^2 + 2XZ - YZ = 0.
\end{align*}
\end{proposition}

\begin{proof}
The only possible decomposition into pairs is $2 \cdot \a{5} + \t$. We assume that this configuration of three conics $Q_1$, $Q_2$ and $Q_3$ is represented by the following graph:
    
    \begin{center}
    \begin{tikzpicture}
    \draw [style={decorate, decoration=snake}] (0,0) -- (2,0);
    \draw (0,0) -- (2,0);
    \draw [style={decorate, decoration=snake}] (0,0) -- (1,1.5);
    \draw (0,0) -- (1,1.5);
    \draw (1,1.5) -- (2,0);
    \draw [fill] (0,0) circle (3pt);
    \node at (-0.5,0) {$Q_1$};
    \draw [fill] (2,0) circle (3pt);
    \node at (2.5,0) {$Q_2$};
    \draw [fill] (1,1.5) circle (3pt);
    \node at (1,2) {$Q_3$};
    \end{tikzpicture}
    \end{center}
We can mimic the proof of Proposition \ref{prop: base}. We use \cite[Proposition 4.2.1]{PHD} to get the equations of $Q_1$ and $Q_2$ along with their parametrizations. Then we write the equation for $Q_3$ in two ways, using Fact \ref{fact:all}~\ref{fact:all_a5} and \ref{fact:all_t}. Then, we compare the coefficients in both equations and as a result we obtain the following formulas:
\begin{align}\label{eq3_param}
\begin{split}
    &Q_1: X^2 - YZ = 0, \\
    &Q_2: X^2 - \left(\frac{1}{2}c^2 + cp\right)Y^2 + cXY - YZ = 0, \\
    &Q_3: (5c + 8p)X^2 - 3(c + 2p)^2XY + 4\left(\frac{1}{2}c + p\right)^3Y^2 - 4XZ + (c + 4p)YZ = 0,
\end{split}
\end{align}
with $c, p \in \mathbb{C}$, $c \neq 0$. This arrangement has the following singularities:
\begin{itemize}
    \item $[0 : 0 : 1]$ -- a $D_8$ singularity, where $Q_1$ and $Q_2$ are tangent,
    \vspace{-6pt}
    \item $\big[\frac{1}{2}c + p : 1 : (\frac{1}{2}c + p)^2\big]$ -- a $D_8$ singularity, where $Q_1$ and $Q_3$ are tangent,
    \vspace{-6pt}
    \item $\big[p : 1 : p^2 - \frac{1}{2}c^2\big]$ -- an $A_3$ singularity, intersection of $Q_2$ and $Q_3$.
\end{itemize}   
In order to use the same method as before, we need a fourth point, uniquely determined by the equations of $Q_1$, $Q_2$ and $Q_3$. We take a point which is an intersection of $Q_3$ and the line tangent to $Q_1$ at $[0 : 0 : 1]$ (different from $[0 : 0 : 1]$). It has coordinates $[4 : 0 : 5c+8p]$. As before, we denote by $M$ the mapping of these $4$ points into $[1:0:0]$, $[0:1:0]$, $[0:0:1]$ and $[1:1:1]$, respectively. Then
$$M^{-1} = \begin{bmatrix} 
0 & 4c+8p & -8p \\ 
0 & 8 & -8 \\ 
-c^2 & 2(c+2p)^2 & -8p^2+4c^2 
\end{bmatrix}.$$
The conics are mapped to
\begin{align*}
    &Q_1: 4Z^2 + XY - XZ - 2YZ = 0, \\
    &Q_2: XY - XZ + 2YZ = 0, \\
    &Q_3: XY + XZ - 2YZ = 0,
\end{align*}
hence all arrangements given by \eqref{eq3_param} are projectively equivalent. We obtain the equations from the statement of our proposition by putting $c=-2$ and $p=1$.

\end{proof}

\begin{figure}[ht]
\centering
\includegraphics[scale=0.6]{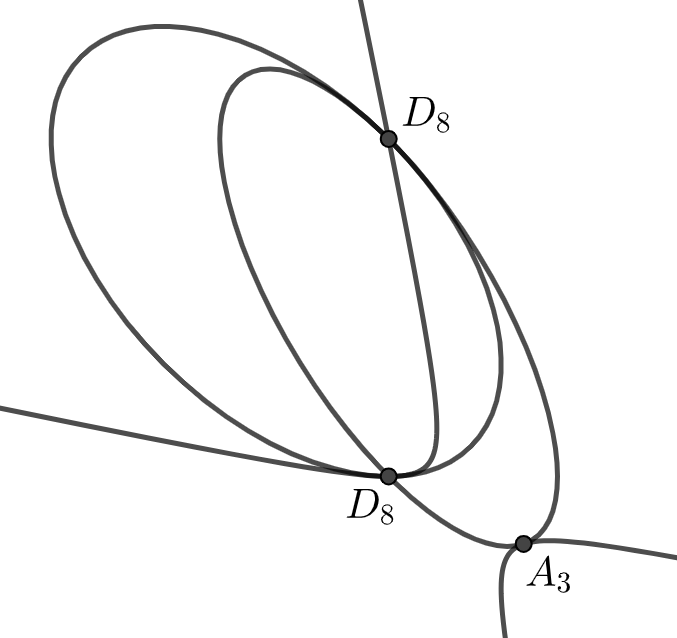}
\caption{An example of the arrangement from Proposition \ref{prop: a3d8d8} \\ ($c = 4$, $p = -2$,\ $Z = 1 - X - Y$).}
\end{figure}

\begin{proposition}
Any free configuration of three conics $Q_1$, $Q_2$, $Q_3$ with weak combinatorics
$$(n_{2}, t_{3}, n_{3}, t_{5}, d_{6}, t_{7}, d_{8}, d_{10}) = (2, 1, 0, 0, 0, 2, 0, 0)$$
is projectively equivalent to the conics
\begin{align*}
    &Q_1: X^2 - YZ = 0, \\
    &Q_2: X^2 - 2Z^2 - YZ = 0, \\
    &Q_3: 2X^2 - Y^2 - 2YZ= 0.
\end{align*}
\end{proposition}

\begin{proof}
The only possible decomposition into pairs is $2\cdot \a{7} + \t$. We assume that this configuration of three conics $Q_1$, $Q_2$ and $Q_3$ is represented by the following graph:

    \begin{center}
    \begin{tikzpicture}
    \draw [line width = 2.5pt] (0,0) -- (2,0);
    \draw [line width = 2.5pt] (0,0) -- (1,1.5);
    \draw (1,1.5) -- (2,0);
    \draw [fill] (0,0) circle (3pt);
    \node at (-0.5,0) {$Q_1$};
    \draw [fill] (2,0) circle (3pt);
    \node at (2.5,0) {$Q_2$};
    \draw [fill] (1,1.5) circle (3pt);
    \node at (1,2) {$Q_3$};
    \end{tikzpicture}
    \end{center}
We can mimic the proof of Proposition \ref{prop: base}. We use \cite[Proposition 4.1.1]{PHD} to get the equations of $Q_1$ and $Q_2$ along with their parametrizations. Then we write the equation for $Q_3$ in two ways, using Fact \ref{fact:all}~\ref{fact:all_t} and \ref{fact:all_a7}. Then, we compare the coefficients in both equations and as a result we obtain the following formulas:
\begin{align}\label{eq4_param}
\begin{split}
    &Q_1: X^2 - YZ = 0, \\
    &Q_2: X^2 - \frac{1}{2}(p-r)^2Z^2 - YZ = 0, \\
    &Q_3: (7p^2 + 2pr - r^2)X^2 - 8pXY + 2Y^2 - 8p^3XZ + (5p^2 - 2pr + r^2)YZ + 2p^4Z^2 = 0,
\end{split}
\end{align}
with $p, r \in \mathbb{C}$, $p \neq r$. This arrangement has the following singularities:
\begin{itemize}
    \item $[0 : 1 : 0]$ -- an $A_7$ singularity, intersection of $Q_1$ and $Q_2$,
    \vspace{-6pt}
    \item $[p : p^2 : 1]$ -- an $A_7$ singularity, intersection of $Q_1$ and $Q_3$,
    \vspace{-6pt}
    \item $\big[r : r^2 - \frac{1}{2}(p-r)^2 : 1\big]$ -- an $A_1$ singularity, intersection of $Q_2$ and $Q_3$,
    \vspace{-6pt}
    \item $\big[2p - r : \frac{1}{2}(7p^2 - 6pr + r^2) : 1\big]$ -- an $A_1$ singularity, intersection of $Q_2$ and $Q_3$,
    \vspace{-6pt}
    \item $\big[p : p^2 - \frac{1}{2}(p-r)^2 : 1]$ -- an $A_3$ singularity, intersection of $Q_2$ and $Q_3$.
\end{itemize}
We denote by $M$ the mapping of the first four of these points (excluding the $A_3$ singularity, because it is collinear with both $A_7$ singularities) into $[1:0:0]$, $[0:1:0]$, $[0:0:1]$ and $[1:1:1]$, respectively. Then
$$M^{-1} = \begin{bmatrix} 
0 & 4p & -2r \\ 
2(p-r)^2 & 4p^2 & p^2-2pr-r^2 \\ 
0 & 4 & -2 
\end{bmatrix}.$$
The conics are mapped to
\begin{align*}
    &Q_1: Z^2 - 4XY + 2XZ + 2YZ = 0, \\
    &Q_2: 2Y^2 + 2XY - XZ - 3YZ = 0, \\
    &Q_3: 2X^2 + 2XY - 3XZ - YZ = 0.
\end{align*}
hence all arrangements given by \eqref{eq4_param} are projectively equivalent. We obtain the equations from the statement of our proposition by putting $p=0$ and $r=2$.
\end{proof}

\begin{remark}
    In this configuration, two $A_7$ singularities and one $A_3$ singularity are collinear.   
\end{remark}

\begin{remark}
    The result presented in the above proposition is consistent with the result shown in \cite[Remark 2.5]{Persson1982}.
\end{remark}

\begin{proposition}
Any free configuration of three conics $Q_1$, $Q_2$, $Q_3$ with weak combinatorics
$$(n_{2}, t_{3}, n_{3}, t_{5}, d_{6}, t_{7}, d_{8}, d_{10}) = (1, 0, 0, 2, 0, 0, 1, 0)$$
is projectively equivalent to the conics
\begin{align*}
    &Q_1: X^2 - YZ = 0, \\
    &Q_2: X^2 + 18(i\sqrt{3} + 3)Y^2 - 18XY - YZ = 0, \\
    &Q_3: 12X^2 + (i\sqrt{3} + 3)XZ - 12YZ = 0.
\end{align*}
\end{proposition}

\begin{proof}
    The only possible decomposition into pairs is $3 \cdot \a{5}$. One $A_5$ singularity coincides with nodes from two other pairs and becomes the $D_8$ singularity. 
    
We can mimic the proof of Proposition \ref{prop: base}. We use \cite[Proposition 4.2.1]{PHD} to get the equations of $Q_1$ and $Q_2$ along with their parametrizations. We assume that these two conics are tangent at a $D_8$ singularity. Then we write the equation for $Q_3$ in two ways, using Fact \ref{fact:all}~\ref{fact:all_a5}. Then, we compare the coefficients in both equations and as a result we obtain the following formulas:
\begin{align}\label{eq5_param}
    \begin{split}
    &Q_1: X^2 - YZ = 0, \\
    &Q_2: X^2 - 3(p-q)(p+\mu)Y^2 + 3(p-q)XY - YZ = 0, \\
    &Q_3: (3p-q+\mu)X^2 - 3p^2XY + p^3Y^2 - XZ + (q-\mu)YZ = 0,
    \end{split}
\end{align}
with $p, q, \mu \in \mathbb{C}$, $p \neq q$ and additionally $(p-q)^2 + 3\mu(p-q+\mu) = 0$. This arrangement has the following singularities:
    \begin{itemize}
        \item $[0 : 0 : 1]$ -- a $D_8$ singularity, where $Q_1$ and $Q_2$ are tangent,
        \vspace{-6pt}
        \item $[p : 1 : p^2]$ -- an $A_5$ singularity, intersection of $Q_1$ and $Q_3$,
        \vspace{-6pt}
        \item $\big[p + \mu : 1 : (p+\mu)^2\big]$ -- an $A_1$ singularity, intersection of $Q_1$ and $Q_2$.
        \vspace{-6pt}
        \item $\big[q : 1 : q^2 - 3(p-q)(p-q+\mu)\big]$ -- an $A_5$ singularity, intersection of $Q_2$ and $Q_3$.
    \end{itemize}
These points can be mapped into $[1:0:0]$, $[0:1:0]$, $[0:0:1]$ and $[1:1:1]$, respectively. If we denote this map by $M$, then
$$M^{-1} = \begin{bmatrix}
0 & (3u - 3)p & (9 - 3u)(p + \mu) \\
0 & 3u - 3 & 9 - 3u \\
(2u - 6)(p-q)^2 & (3u - 3)p^2 & (9 - 3u)(p + \mu)^2
\end{bmatrix},$$
with $\mu = \frac{u + 3}{6}(q - p)$, either $u = i\sqrt{3}$ or $u = -i\sqrt{3}$. As a result, we obtain two arrangements, depending on the value of $u$. Similarly to the proof of Proposition \ref{prop:twocases}, it can be shown that these arrangements are projectively equivalent, hence all arrangements given by \eqref{eq5_param} are projectively equivalent. We obtain the equations from the statement of our proposition by putting $p = 0$, $q = 6$ and $\mu = i\sqrt{3} + 3$ in \eqref{eq5_param}.
\end{proof}

\begin{remark}
It is easy to verify that the equality $(p-q)^2 + 3\mu(p - q + \mu) = 0$ cannot be satisfied for $p, q, \mu \in \mathbb{R}$ (and $p \neq q$). Therefore this arrangement cannot be realized over $\mathbb{R}$.    
\end{remark}

\begin{proposition}\label{prop: last}
Any free configuration of three conics $Q_1$, $Q_2$, $Q_3$ with weak combinatorics
$$(n_{2}, t_{3}, n_{3}, t_{5}, d_{6}, t_{7}, d_{8}, d_{10}) = (1, 1, 0,1, 0, 0, 0, 1)$$
is projectively equivalent to the conics
\begin{align*}
    &Q_1: X^2 - YZ = 0, \\
    &Q_2: X^2 - 3Z^2 - YZ = 0, \\
    &Q_3: 2X^2 + 3XY - 2YZ = 0.
\end{align*}
\end{proposition}

\begin{proof}
    The only possible decomposition into pairs is $\a{7} + \a{5} + \t$. We assume that this configuration of three conics $Q_1$, $Q_2$ and $Q_3$ is represented by the following graph:
    
    \begin{center}
    \begin{tikzpicture}
    \draw [line width = 2.5pt] (0,0) -- (2,0);
    \draw (0,0) -- (1,1.5);
    \draw [style={decorate, decoration=snake}] (0,0) -- (1,1.5);
    \draw (1,1.5) -- (2,0);
    \draw [fill] (0,0) circle (3pt);
    \node at (-0.5,0) {$Q_1$};
    \draw [fill] (2,0) circle (3pt);
    \node at (2.5,0) {$Q_2$};
    \draw [fill] (1,1.5) circle (3pt);
    \node at (1,2) {$Q_3$};
    \end{tikzpicture}
    \end{center}

We can mimic the proof of Proposition \ref{prop: base}. We use \cite[Proposition 4.1.1]{PHD} to get the equations of $Q_1$ and $Q_2$ along with their parametrizations. Then we write the equation for $Q_3$ in two ways, using Fact \ref{fact:all}~\ref{fact:all_a5} and \ref{fact:all_t}. Then, we compare the coefficients in both equations and as a result we obtain the following formulas:
\begin{align}\label{eq6_param}
    \begin{split}
    &Q_1: X^2 - YZ = 0, \\
    &Q_2: X^2 - 3(p-q)^2Z^2 - YZ = 0, \\
    &Q_3: (-8p+2q)X^2 + 3XY + 9p^2XZ + (-p-2q)YZ - 3p^3Z^2 = 0,
    \end{split}
\end{align}
with $p, q \in \mathbb{C}$ and $p \neq q$. This arrangement has the following singularities:
    \begin{itemize}
        \item $[0 : 1 : 0]$ -- a $D_{10}$ singularity, where $Q_1$ and $Q_2$ are tangent,
        \vspace{-6pt}
        \item $[p : p^2 : 1]$ -- an $A_5$ singularity, intersection of $Q_1$ and $Q_3$,
        \vspace{-6pt}
        \item $[q : q^2 - 3(p-q)^2 : 1]$ -- an $A_3$ singularity, intersection of $Q_2$ and $Q_3$,
        \vspace{-6pt}
        \item $[3p-2q : 6p^2 - 6pq + q^2 : 1]$ -- an $A_1$ singularity, intersection of $Q_2$ and $Q_3$.
    \end{itemize}
These points can be mapped into $[1:0:0]$, $[0:1:0]$, $[0:0:1]$ and $[1:1:1]$, respectively. If we denote this map by $M$, then
$$M^{-1} = \begin{bmatrix}
0 & 3p & -2q \\
-3(p-q)^2 & 3p^2 & 6p^2-12pq+4q^2 \\
0 & 3 & -2
\end{bmatrix}.$$
As a result, we obtain the equations not depending on $p$ and $q$, therefore all arrangements given by \eqref{eq6_param} are projectively equivalent. We obtain the equations from the statement of our proposition by putting $p = 0$ and $q = 1$.
\end{proof}

\begin{figure}[ht!]
\centering
\includegraphics[scale=0.6]{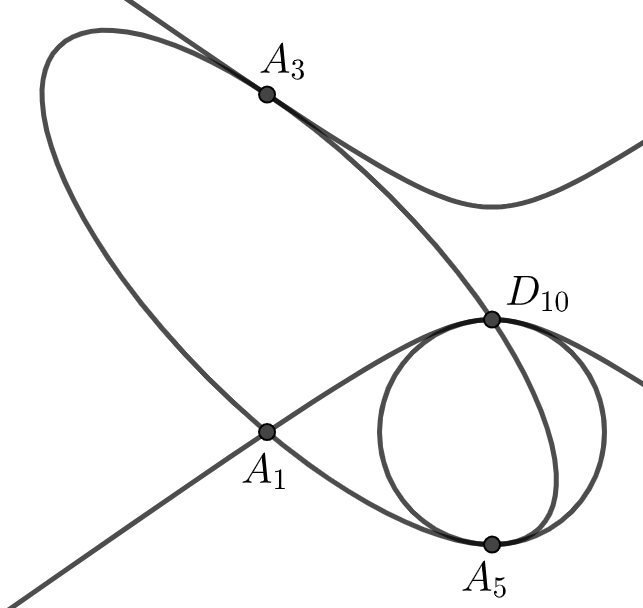}
\caption{An example of the arrangement from Proposition \ref{prop: last} \\ ($p = 0,\ q = 1,\ Z = 1 - Y$).}
\end{figure}

It can be checked that in the above propositions, from Proposition \ref{prop: first} to Proposition \ref{prop: last}, we analyzed all possible weak combinatorics listed in Theorem \ref{thm: weakcombADE}. Therefore we can gather all these results and write them down in a single main theorem:

\begin{thmABC}\label{thm:main1}
Every free arrangement of three smooth conics $Q_1, Q_2, Q_3$ over $\mathbb{C}$ with only ${\rm ADE}$ singularities is projectively equivalent to one of the following $6$ arrangements:

\begin{enumerate}

    \item [(i)]
        \begin{enumerate}
            \item[] $Q_1: X^2 - YZ = 0,$
            \item[] $Q_2: X^2 - 2XY - YZ = 0,$
            \item[] $Q_3: X^2 + 2XZ - YZ = 0.$
        \end{enumerate}
        This arrangement has one $A_3$ singularity and two $D_8$ singularities.

    \item [(ii)]       
    \begin{enumerate}
            \item [] $Q_1: X^2 - YZ = 0$, 
            \item [] $Q_2: 3X^2 + Y^2 + 3XY - 3YZ = 0$, 
            \item [] $Q_3: X^2 + 27Z^2 + 9XZ - YZ = 0.$
        \end{enumerate}
        This arrangement has one $D_4$ singularity and three $A_5$ singularities.

    \item[(iii)]
        \begin{enumerate}
            \item[] $Q_1: X^2 - YZ = 0,$
            \item[] $Q_2: X^2 + i\sqrt{3}Z^2 - YZ = 0,$
            \item[] $Q_3:8X^2 - 3(i\sqrt{3}+1)Y^2 + 8i\sqrt{3}XY - 8YZ = 0.$
        \end{enumerate}
        This arrangement has two $A_1$, two $A_5$ and one $A_7$ singularities. 

    \item[(iv)]
        \begin{enumerate}
            \item[] $Q_1: X^2 - YZ = 0,$
            \item[] $Q_2: X^2 - 2Z^2 - YZ = 0,$
            \item[] $Q_3 : 2X^2 - Y^2 - 2YZ= 0.$
        \end{enumerate}
        This arrangement has two $A_1$ singularities, one $A_3$ singularity and two $A_7$ singularities.

    \item[(v)]
    \begin{enumerate}
        \item[] $Q_1 : X^2 - YZ = 0,$
        \item[] $Q_2 : X^2 + 18(i\sqrt{3}+3)Y^2 - 18XY - YZ = 0,$
        \item[] $Q_3 : 12X^2 + (i\sqrt{3}+3)XZ - 12YZ = 0.$
    \end{enumerate}
    This arrangement has one $A_1$, two $A_5$ and one $D_8$ singularity.

    \item[(vi)]
    \begin{enumerate}
        \item[] $Q_1 :  X^2 - YZ = 0,$
        \item[] $Q_2 : X^2 - 3Z^2 - YZ = 0,$
        \item[] $Q_3 : 2X^2 + 3XY - 2YZ = 0.$
    \end{enumerate}
    This arrangement has one $A_1$, one  $A_3$, one $ A_5$ and one $D_{10}$ singularity. 
\end{enumerate}
\end{thmABC}

\section{Adding the $J_{2,0}$ singularity}

Given the complete classification of free arrangements of three smooth conics with only ${\rm ADE}$ singularities, a natural question is to ask how the classification looks when considering all possible quasi-homogeneous singularities. The first step towards the full classification is considering another type of quasi-homogeneous singularity, namely $J_{2,0}$. For this singularity, there are three conics passing through one point such that each pair of conics intersects as an $A_3$ singularity. We present below an appendix to Table \ref{table: ADES} with the information about the $J_{2,0}$ singularity.

\begin{table}[ht]
    \centering
    \begin{tabular}{c|c|c|c}
     & Number of & & \\
     Singularity & occurrences & $\mu_{p}$ & Multiplicity \\ \hline
    $J_{2,0}$ & $j$ & 10 & 6
    \end{tabular}
    \caption {The $J_{2,0}$ singularity.}
    \label{table: QHS}
\end{table}

By repeating the reasoning from the previous section, we obtain the following system of equations, which is a necessary condition for such free arrangements to exist:
\begin{equation}\label{system_equations_2}
\begin{cases}
    n_{2} + 3t_{3} + 4n_{3} + 5t_{5} + 6d_{6} + 7t_{7} + 8d_{8} + 10d_{10} + 10j = \tau(\mathcal{C}) = 25 - 5d_1 + d_1^2,\\
    n_{2} + 2t_{3} + 3n_{3} + 3t_{5} + 4d_{6} +  4t_{7} + 5d_{8} + 6d_{10} + 6j = 12,
\end{cases}
\end{equation}
Here, the right side of the first equation comes from Theorem \ref{duPles} for $d = 6$, where the second equation refers to the well-known naive combinatorial counting of the intersections of conics. We have $d_1 \in \{1, 2\}$. It turns out that \eqref{system_equations_2} has exactly one solution when $d_1 = 1$, namely 
$$(n_{2}, t_{3}, n_{3}, t_{5}, d_{6}, t_{7}, d_{8}, d_{10}, j) = (0, 0, 0, 0, 0, 3, 0, 0, 0),$$
however, it cannot be realized over $\mathbb{C}$ by smooth conics (see \cite[Corollary 1.2]{Zielinski2023}). We also obtain $6$ new solutions (apart from the ones from Theorem \ref{thm: weakcombADE}) when $d_1 = 2$, so we can write the following result:

\begin{proposition}\label{prop: weakcombAll} A configuration of three smooth conics with only ${\rm ADE}$ singularities and at least one $J_{2,0}$ singularity is free if and only if it has one of the following weak combinatorics:
 {\small
\begin{align*}
(n_{2}, t_{3}, n_{3}, t_{5}, d_{6}, t_{7}, d_{8}, d_{10}, j) \in \{(0, 0, 1, 1, 0, 0, 0, 0, 1), (0, 1, 0, 0, 1, 0, 0, 0, 1), \\
(0, 3, 0, 0, 0, 0, 0, 0, 1), (1, 0, 0, 0, 0, 0, 1, 0, 1), (1, 1, 0, 1, 0, 0, 0, 0, 1), (2, 0, 0, 0, 0, 1, 0, 0, 1)\}.
\end{align*}
}
\end{proposition}

In the following, we will show that none of the above weak combinatorics is realizable over $\mathbb{C}$. To begin with, let us write down a useful lemma.

\begin{lemma} \label{split_noj}
    If a configuration of three smooth conics has at least one singularity of type $A_5$, $A_7$, $D_8$ or $D_{10}$, then it does not contain a singularity of type $J_{2,0}$.
\end{lemma}

\begin{proof}
Suppose that a configuration of three smooth conics ${Q_1, Q_2, Q_3}$ has a $D_8$ singularity at a point $P_1$ and a $J_{2,0}$ singularity at a point $P_2$. Without loss of generality, assume that $Q_1$ and $Q_2$ intersect at $P_1$ with multiplicity $3$. Then these two conics must also intersect at $P_2$ with multiplicity $2$, which contradicts B\'ezout's Theorem. The proof in the remaining cases ($A_5$, $A_7$, and $D_{10}$) is analogous.  
\end{proof}

\begin{proposition}
    The following weak combinatorics
\begin{align*}
(n_{2}, t_{3}, n_{3}, t_{5}, d_{6}, t_{7}, d_{8}, d_{10}, j) \in \{(0, 0, 1, 1, 0, 0, 0, 0, 1), \\
(1, 0, 0, 0, 0, 0, 1, 0, 1), (1, 1, 0, 1, 0, 0, 0, 0, 1), (2, 0, 0, 0, 0, 1, 0, 0, 1)\}.
\end{align*}
    cannot be realized over $\mathbb{C}$ by smooth conics.
\end{proposition}
\begin{proof}
    Follows from Lemma \ref{split_noj}.   
\end{proof}

\begin{proposition}
    The following weak combinatorics
    $$(n_{2}, t_{3}, n_{3}, t_{5}, d_{6}, t_{7}, d_{8}, d_{10}, j) = (0, 1, 0, 0, 1, 0, 0, 0, 1)$$
    cannot be realized over $\mathbb{C}$ by smooth conics.
\end{proposition}
\begin{proof}
    Follows from Lemma \ref{split_nod2}.
\end{proof}

\begin{proposition}
    The following weak combinatorics
    $$(n_{2}, t_{3}, n_{3}, t_{5}, d_{6}, t_{7}, d_{8}, d_{10}, j) = (0,3,0,0,0,0,0,0,1)$$
    cannot be realized over $\mathbb{C}$ by smooth conics.
\end{proposition}

\begin{proof}
    Suppose that a configuration realizing this weak combinatorics exists. The only possible decomposition into pairs is $3 \cdot \tt$. By \cite[Proposition 3]{Megyesi2000}, we may assume that
$$Q_1: X^2 + Y^2 - Z^2 = 0, \quad Q_2: \frac{1}{q^2}X^2 + Y^2 - Z^2 = 0,$$
for some $q \in \mathbb{C}$, $q \notin \{0, 1, -1\}$. These conics intersect in two points: $[0 : 1 : 1]$ and $[0 : -1 : 1]$. Without loss of generality we may assume that $[0 : 1 : 1]$ is the $J_{2,0}$ singularity.
By \cite[Remark 4.3.2]{PHD} we know that $Q_1$, $Q_2$ have the following parametrizations:
\begin{align}
\label{param1} Q_1 &= \{[2st : t^2 - s^2 : t^2 + s^2 ] \mid [s:t] \in \mathbb{P}^1 \}, \\
\label{param2} Q_2 &= \{[2qst  : t^2-s^2:t^2 + s^2] \mid  [s:t] \in \mathbb{P}^1 \}.
\end{align}
By Fact \ref{fact:all}~\ref{fact:all_t2} we know that $Q_3$ is given by the equation $\lambda Q_1 + \ell_1 \ell_2 = 0$, where  $\lambda \neq 0$ and $\ell_1,\ell_2$ are lines tangent to $Q_1$ at two intersection points. We have $[0:1:1] \in Q_1 \cap Q_3$ and from \eqref{param1} the second point of intersection is of the form $[2u : u^2 - 1 : u^2 + 1]$ for some $u \in \mathbb{C}$, $u \neq 0$. We therefore have $\ell_1 = Y- Z$ and $ \ell_2 = 2uX + (u^2 -1)Y - (u^2 +1)Z$. In the end we write
$$
Q_3 : \lambda(X^2 + Y^2 -Z^2) + (Y-Z)(2uX + (u^2 - 1)Y - (u^2 + 1)Z ) = 0.
$$
Now, consider the intersection of $Q_2$ and $Q_3$. From \eqref{param2}, we can write
$$Q_2 \setminus \{[0 : -1 : 1]\} = \{[2qv : 1 - v^2 : 1+ v^2] \mid v \neq 0 \}.$$
Let $f_{23}$ denote the polynomial in variable $v$ obtained by substituting the above parametrization to the above equation of $Q_3$. We have
\begin{equation}\label{eq:f23}
f_{23}(v) = 4v^2(u^2v^2 - 2uqv + (\lambda q^2 - \lambda + 1)).
\end{equation}
Since $[0:-1:1]$ cannot lie in the intersection of $Q_2$ and $Q_3$, the roots of $f_{23}$ are simply values of $v$ that give all the intersection points of $Q_2$ and $Q_3$. From the intersection behavior of $Q_2, Q_3$ we also know that $f_{23}$ must be of the form
$$f_{23}(v) = A(v-\mu_1)^2(v - \mu_2)^2.$$
for $\mu_1, \mu_2 \in \mathbb{C}$, $\mu_1 \neq \mu_2$. We know that the $J_{2,0}$ singularity lies in the in the intersection of $Q_2$ and $Q_3$ therefore we know that $\mu_1 = 0$, which leads to
$$f_{23}(v) = Av^2(v-\mu_2)^2.$$
Thus, by comparing the above equation to \eqref{eq:f23}, we get that the discriminant of $u^2v^2 - 2uqv + (\lambda q^2 - \lambda + 1)$ is equal to $0$. Therefore we can write
$$4 u^2 q^2 - 4 u^2(\lambda q^2 - \lambda + 1) =  -4u^2 (\lambda - 1)(q^2 -1) = 0$$
and hence $\lambda = 1$ or $u = 0$ or $q \in \{-1,1\}$. Since $u = 0$, $q \in \{ -1,1 \}$ are all contradictory with our assumptions, the only valid option is $\lambda = 1$. If we substitute $\lambda = 1$ into the equation of $Q_3$ we get
$$
X^2 + Y^2 - Z^2 + (Y-Z)(2uX + (u^2 -1)Y - (u^2-1)Z) = (X + uY - uZ)^2,
$$
which is a contradiction with irreducibility of $Q_3$. \end{proof}

As a result from the above propositions, we obtain the second main theorem:
\begin{thmABC}\label{thm:main2}
Every free arrangement of three smooth conics $Q_1, Q_2, Q_3$ over $\mathbb{C}$ admitting only ${\rm ADE}$ singularities and $J_{2,0}$ singularities is projectively equivalent to one of the $6$ arrangements presented in Theorem \ref{thm:main1}.
\end{thmABC}

\bibliographystyle{abbrv}
\bibliography{master}

\vspace{10pt}
Affiliation of the authors: Department of Mathematics, University of the National Education Commission, Krakow, Podchor\k{a}\.{z}ych 2, PL-30-084 Krak\'{o}w, Poland.

\L{}ukasz Merta: \texttt{lukasz.merta@uken.krakow.pl}

Filip Zieli\'nski: \texttt{filipzielinski377@gmail.com}

Marcin Zieli\'nski: \texttt{marcin.zielinski@uken.krakow.pl}

\end{document}